\documentclass{article}
\usepackage{amssymb}
\usepackage{latexsym}

\input{tcilatex}
\begin{document}

\begin{center}
\textbf{\LARGE From Infinitesimal Harmonic}

\textbf{\LARGE Transformations to Ricci Solitons}

\textbf{\Large Sergey E. Stepanov}
\end{center}

\bigskip

\begin{center}
\textit{\small Department of Mathematics, Finance University under the
Government of Russian Federation,} \textit{49, Leningradsky Prospect,
Moscow, 125993 Russia}

\textbf{{\_\_}{\_}{\_}{\_}{\_}{\_}{\_}{\_}{\_}{\_}{\_}{\_}{\_}{\_}{\_}{\_}{\_%
}{\_}{\_}{\_}{\_}{\_}{\_}{\_}{\_}{\_}{\_}{\_}{\_}{\_}{\_}{\_}{\_}{\_}{\_}{\_}%
{\_}{\_}{\_}{\_}{\_}{\_}{\_}{\_}{\_}{\_}{\_}{\_}{\_}{\_}{\_}{\_}{\_}{\_}{\_}{%
\_}{\_}{\_}{\_}{\_}{\_}{\_}{\_}{\_}{\_}{\_}{\_}}
\end{center}

\bigskip \textbf{Abstract}

{\small The concept of the Ricci soliton was introduced by Hamilton. Ricci
soliton is defined by vector field and it's a natural generalization of
Einstein metric. We have shown earlier that the vector field of Ricci
soliton is an infinitesimal harmonic transformation. In our paper, we survey
Ricci solitons geometry as an application of the theory of infinitesimal
harmonic transformations.} \bigskip

 \bigskip


\begin{center}
\textbf{1. Harmonic diffeomorphisms and infinitesimal harmonic
transformations}
\end{center}

A smooth mapping $f:\left( {M, g} \right)\to \left( {{M}^{\prime }, {g}%
^{\prime }} \right)$ between two Riemannian manifolds is called \textit{%
harmonic} (see [3]) if $f$ provides an exstremum of the Dirichlet functional 
$E_\Omega \left( f \right)=\frac{1}{2}\int_{_\Omega } {\left\| {\mbox{d}f}
\right\|} ^2\mbox{d}V$ with respect to the variations of $f$ that are
compactly supported in a relatively compact open subset $\Omega \subset M$.
(Here, d$V $ is the volume element of the metric $g$.) The following theorem
is true (see [3]).

\textbf{Theorem 1.1}. A smooth mapping $f:\left( {M, g} \right)\to \left( {{M%
}^{\prime }, {g}^{\prime }} \right)$ is harmonic if and only if it satisfies
the Euler-Lagrange equations 
$$
\label{eq1} g^{ij\,}\left( {\partial_i \partial_j f^\beta -\Gamma_{ij}^k
\partial_k f^\beta +\partial_i f^\beta \partial_j f^\gamma \left( {{%
\Gamma^{\prime }}_{\beta \gamma }^\alpha \circ f} \right)} \right)=0\eqno%
(1.1) 
$$
where $y^\alpha =f^\alpha \left( {x^1,...,x^n} \right)$ is local
representation of $f$; $g^{ij}$ are local contrvariant components of the
metric tensor $g$; $\Gamma_{ij}^k $ and ${\Gamma^{\prime }}_{\beta \gamma
}^\alpha$ are Christoffel symbols of $( {M, g} )$ and $( {{M}^{\prime }, {g}%
^{\prime }} )$ respectively; $i,j,k=1,\dots,n=\mathrm{dim}M$ and $%
\alpha,\beta,\gamma=1,\dots,n^{\prime }=\mathrm{dim}M^{\prime }$.

If we suppose that $\mathrm{dim}M=\mathrm{dim}M^{\prime }=n$ and $f:\left( {%
M, g} \right)\to \left( {{M}^{\prime }, {g}^{\prime }} \right)$ is a
diffeomorphism then $f $ is locally represented by the following equations $%
y^{i} = x^{i\thinspace \thinspace }$for $i$, \textit{j, k ,{\ldots}} = 1, 2, 
{\ldots} , $n $ and therefore the Euler-Lagrange equations (1.1) take the
form

$$
g^{ij}((\Gamma^{\prime k}_{ij}\circ f)-\Gamma^{k}_{ij})=0\eqno(1.2) 
$$
where $\Gamma^{k}_{ij}$ and $\Gamma^{\prime k}_{ij}$ are the Christoffel
symbols of the Levi-Civita connection $\nabla$ on $( {M, g} )$ and $%
\nabla^{\prime }$ on $( {M^{\prime }, g^{\prime }} )$ respectively.

{Suppose that we have a local one-parameter group of infinitesimal point
transformations $f_{t}(x)=x^{\prime k}+t\xi ^{k})$ generated by a vector
field $\xi =\xi ^{k}\partial _{k}$ on $({M,g})$ for so-called canonical
parameter t such that $t\in (-\xi ,+\xi )\subset \mathbf{R}$. In this case
the Lie derivative of the Christoffel symbols $\Gamma _{ij}^{k}$ of the
Levi-Civita connection $\nabla $ has the form (see [21], pp. 8-9)} 
$$
\label{eq2}(L_{\xi }\Gamma _{ij}^{k})t=\Gamma _{ij}^{\prime k}-\Gamma
_{ij}^{k}=\nabla _{i}\nabla _{j}\xi ^{k}-R_{ijl}^{k}\xi ^{l}\eqno(1.3) 
$$%
where $\Gamma _{ij}^{\prime k}(x)=f_{t}^{\ast }(\Gamma _{ij}^{k}(x^{\prime
}))$.

\textbf{Definition} (see [14]; [19]). A vector field $\xi $ on $(M,g)$ is
called an \textit{infinitesimal harmonic transformation }if the
one-parameter group of local transformations of $(M,g)$ generated by $\xi $
consists of local harmonic diffeomorphisms.

{By the definition and (1.3) we deduce the following equation } 
$$
\label{eq3} \Delta\theta=2Ric^*\xi\eqno(1.4) 
$$
where $\xi$ is an infinitesimal harmonic transformation and $%
\theta=g(\xi,\cdot)$ is its dual 1-form; $\Delta:=\mathrm{dd}^*+\mathrm{d}^*%
\mathrm{d}$ is the Hodge Laplacian on the space 1-forms $\Omega^1(M)$; $%
\mathrm{Ric}^*$ is the linear Ricci operator defined by the identity $g(%
\mathrm{Ric}^*X,\cdot)=\mathrm{Ric}(X,\cdot)$ for the tensor Ricci $\mathrm{%
Ric}$ and an arbitrary vector field $Y$ on $M$.

\textbf{Theorem 1.2} (see [14], [19]). \textit{The equation} $\Delta \theta
=2\mathrm{Ric}^{\ast }\xi $ \textit{is a necessary and sufficient condition
for vector field }$\xi $\textit{\ to be an infinitesimal harmonic
transformation on a Riemannian manifold }$(M,g).$

\begin{center}
\textbf{2. Examples of infinitesimal harmonic transformations}
\end{center}

In this paragraph we will give five examples of infinitesimal harmonic
transformations on Riemannian, nearly Kahlerian and Kahlerian manifolds.

\textbf{Example 1}. An infinitesimal isometric transformation on a
Riemannian manifold is an infinitesimal harmonic transformation.

A vector field $\xi $ on an $n$-dimensional Riemannian manifold $(M,g)$ is
an \textit{infinitesimal isometric transformation} if $L_{\xi }g=0$ where $%
L_{\xi }$ is the Lee derivative in direction to $\xi $. By direct
computation, we can deduce the following equalities $\Delta \theta =2\mathrm{%
Ric}^{\ast }\xi $ and $\mathrm{d^{\ast }\theta }=0$ for $\theta =g(\xi
,\cdot )$. Moreover, these equalities are a necessary and sufficient
condition for a vector field $\xi $ to be an infinitesimal isometric
transformation on a compact Riemannian manifold $(M,g)$ (see [21], p. 221).

\textbf{Example 2}. An infinitesimal conformal transformation on a
two-dimensional Riemannian manifold is a harmonic transformation.

Recall that a vector field $\xi $ is an \textit{infinitesimal conformal
transformation} if $L_{\xi }g=-\frac{2}{n}(\mathrm{d}^{\ast }\theta )g$ for $%
\theta =g(\xi ,\cdot )$. By direct computation, we can deduce the following
equality $\Delta \theta +(1-\frac{2}{n})\mathrm{dd}^{\ast }\theta =2\mathrm{%
Ric}^{\ast }\xi $. Moreover, by virtue of the Lihnerowicz theorem (see [12])
this equality is a necessary and sufficient condition for a vector field $%
\xi $ to be an infinitesimal conformal transformation on a compact
Riemannian manifold $(M,g)$. In particular, for $n=2$ we have the equality $%
\Delta \theta =2\mathrm{Ric}^{\ast }\xi $. Therefore, any infinitesimal
harmonic transformation on a two-dimension compact Riemannian manifold is an
infinitesimal conformal transformation.

\textbf{Example 3 }(see [19])\textbf{.} A holomorphic vector field on a
nearly Kahlerian manifold is infinitesimal harmonic transformation.

Let the triplet $(M,g,\mathrm{J})$ be a \textit{nearly Kahlerian manifold}
(see [6]) where $\mathrm{J}\in T^{\ast }M\otimes TM$ such that $\mathrm{J}%
^{2}=-\mathrm{id}_{M}$, $g(\mathrm{J},\mathrm{J})=g$ and $(\nabla _{X}%
\mathrm{J})Y+(\nabla _{Y}\mathrm{J})X=0$ for any $X,Y\in TM$ and let $\xi $
be a \textit{holomorphic vector field }on $(M,g,\mathrm{J})$, i.e. $L_{\xi }%
\mathrm{J}=0$. In this case as we have proved in [19] that the following
identity $\Delta \theta =2\mathrm{Ric}^{\ast }\xi $ holds.

\textbf{Remark 4.} On a compact Kahlerian manifold $(M,g,\mathrm{J})$, where
as well known $\nabla \mathrm{J}=0$, a vector field $\xi $ is holomorphic if
and only if $\Delta \theta =2\mathrm{Ric}^{\ast }\xi $ (see [21], p. 280).
Therefore, in particular, a vector field $\xi $ on a compact Kahlerian
manifold is an infinitesimal harmonic transformation if and only if $\xi $
is holomorphic.

\textbf{Example 5 }(see [20]). A vector field $\xi $ that makes a Riemannian
metric $g$ into a Ricci soliton metric is necessarily an infinitesimal
harmonic transformation.

Let $M$ be a smooth manifold. A \textit{Ricci soliton} $(g,\xi,\lambda)$ is
a Riemannian metric $g$ together with a vector field $\xi$ on $M $ and some
constant $\lambda$ that satisfies the equation $-2\mathrm{Ric} =L_\xi
g+2\lambda g$ (see [1], pp 22-23).

The Lie derivative of $\nabla $ has the following form (see [21], p. 52) 
$$
L_{\xi }\Gamma _{ij}^{k}=\frac{1}{2}g^{kl}(\nabla _{i}L_{\xi }g_{jl}+\nabla
_{j}L_{\xi }g_{il}-\nabla _{l}L_{\xi }g_{ij}).\eqno(2.1) 
$$%
Substituting the identity $L_{\xi }g=-2(\mathrm{Ric}+\lambda g)$ in (3.1) we
find $L_{\xi }\Gamma _{ij}^{k}=g^{kl}(-\nabla _{i}R_{jl}-\nabla
_{j}R_{il}+\nabla _{l}R_{ij})$ for local components $R_{ij}$ of the Ricci
tensor $\mathrm{Ric}$. From the last equation we have $g^{ij}(L_{\xi }\Gamma
_{ij}^{k})=g^{kl}(-2\nabla _{j}R_{l}^{j}+\nabla _{l}s)=0$ for the scalar
curvature $s=g^{ij}R_{ij}$. Here we have taken advantage of the Schur's
lemma $2\nabla _{j}R_{l}^{j}=\nabla _{l}s$.

\textbf{Remark}. If $\theta=\mathrm{d}F$ for a smooth function $%
F:M\rightarrow \mathbf{R}$ then the equation of an infinitesimal harmonic
transformation $\Delta \theta=2\mathrm{Ric}^*\xi$ can be written as $%
\Delta(\nabla_kF)=2R_k^j\nabla_jF$ where $\Delta(\nabla_kF)=\nabla_k(\Delta
F)$. On the other hand, if we put $\xi=\mathrm{grad\;F}$ then from the
equation of a Ricci soliton we conclude $\Delta F=s+n\lambda$ and hence the
equation $\nabla_k(\Delta F)=2R_k^j\nabla_jF$ is equal to $\nabla_k
s=2R_k^j\nabla_j F$. The last equation was proved by Hamilton for a gradient
Ricci soliton (see [8]).

\begin{center}
\textbf{3. The Yano Laplacian}
\end{center}

Let $(M,g)$ be a compact Riemannian manifold. We may also assume that $(M,g)$
is orientable; if $(M,g)$ is not orientable, we have only to take an
orientable twofold covering space of $(M, g)$. Denote by $S^pM$ the bundle
of symmetric bilinear forms on $(M, g)$, $\delta^*$ the symmetric
differentiation operator $\delta^*:C^\infty S^pM\rightarrow C^\infty S^{p+1}
M$ and $\delta$ the linear differential operator $\delta : C^\infty
S^{p+1}M\rightarrow C^\infty S^pM$ as the adjoint operator to $\delta^*$
with respect to the global scalar product on $S^pM$ 
\[
\langle \varphi,\varphi^{\prime }\rangle =\int_M\frac{1}{p!}%
g(\varphi,\varphi^{\prime })dV, 
\]
which we get by integrating the pointwise inner product $g(\varphi,\varphi^{%
\prime })$ for all $\varphi,\varphi^{\prime \infty }S^pM$.

\textbf{Definition} (see [17]; [18]). A differential operator $\Box
:C^{\infty }S^{p}M\rightarrow C^{\infty }S^{p}M$ is called the \textit{Yano
differential operator }if $\Box =\delta \delta ^{\ast }-\delta ^{\ast
}\delta $.

The Yano operator $\Box $ and Bochner Laplacian $\nabla ^{\ast }\nabla $ are
connected by the Weitzenbock formula $\Box =\nabla ^{\ast }\nabla +\Re _{p}$
for the symmetric endomorphism $\Re _{p}$ of the bundle $S^{p}M$ such that $%
\Re _{p}$ can be algebraically (even linearly) expressed through the
curvature and Ricci tensors of $(M,g)$ (see [17]; [18]). In particular for $%
p=1$ we have $\Re _{1}=-\mathrm{Ric}^{\ast }$ and since $\Box =\Delta -2%
\mathrm{Ric}^{\ast }$ (see [18]).

\textbf{Remark}. This form of the operator $\Box $ was used by K. Yano (see
[22], p. 40) for the investigation of local isometric transformations of $%
(M,g)$. Therefore we have named $\Box $ as the Yano operator. Moreover, Yano
has named a vector field $\xi $ as \textit{geodesic} if $\Box \xi =0$ (see
[23]).

In view of told above we can formulate the following theorem.

\textbf{Theorem 3.1}. (see [19]). \textit{A necessary and sufficient
condition for a vector field } $\xi $ \textit{\ on a Riemannian manifold }$%
(M,g)$\textit{\ to be infinitesimal harmonic transformation is that} $\xi
\in \mathrm{Ker}\Box $ \textit{for the Yano operator} $\Box $.

From the following identity $\langle \Box \varphi ,\varphi ^{\prime }\rangle
=\langle \varphi ,\Box \varphi ^{\prime }\rangle $, we conclude that $\Box $
is a self-adjoint differential operator (see [15]). In addition, the symbol 
\textit{$\sigma $ }of the Yano operator $\Box $ satisfies (see [17]) the
following condition \textit{$\sigma (\Box )(\vartheta ,x)\varphi
_{x}=-g(\vartheta ,\vartheta )\varphi _{x}$} for an arbitrary $x\in M$ and $%
\vartheta \in T_{x}^{\ast }M-\{0\}$. Hence the Yano operator $\Box $ is the
self-adjoint Laplacian operator and its kernel is a finite-dimensional
vector space on compact $(M,g)$. In addition we recall that by virtue of the
Fredholm alternative (see [16], p. 205) the vector spaces Ker$\Box $ and Im$%
\Box $ are orthogonal complement of each other with respect the global
scalar product defined on compact $(M,g)$, i.e. $\Omega ^{p}(M)=\mathrm{Ker}%
\Box \oplus \mathrm{Im}\Box $. In particular, for $p=1$ we can formulate the
following

\textbf{Theorem 3.2}. \textit{The vector space} Ker$\Box$ \textit{of all
infinitesimal harmonic transformations on compact Riemannian manifold }$%
(M,g) $\textit{\ is a finite-dimensional vector space and the following
orthogonal decomposition} $\Omega^1(M)= \mathrm{Ker} \Box\oplus \mathrm{Im}%
\Box$ \textit{holds.}

For any conformal Killing vector field $\zeta $ and its dual 1-form $\omega $
on compact smooth manifold $(M,g)$ we have $\langle \Delta \omega
+(1-2n^{-1})\mathrm{dd}^{\ast }\omega -2\mathrm{Ric}^{\ast }\zeta ,\omega
\rangle \geq 0$ (see [9]). From this inequality we conclude that $\langle
\delta ^{\ast }\omega ,\delta ^{\ast }\omega \rangle \geq 2n^{-1}\langle
\delta \omega ,\delta \omega \rangle \geq 0$ and hence $\langle \Box \omega
,\omega \rangle \geq 0$ for $n\geq 2$.

\textbf{Remark}. An infinitesimal harmonic transformation $\xi $\ is a
harmonic vector field if and only if \ $\mathrm{Ric}^{\ast }\xi =0$. On a
compact Riemannian manifold this infinitesimal harmonic transformation $\xi $%
\ must be a covariant constant vector field. Therefore, in particular, if
Ricci tensor is a nosingular tensor then determinant of the $\mathrm{Ricci}$
tensor\ is nonzero for every point $x\in M$ and in this case does not exist
a nonzero harmonic vector field that belongs to the vector space $\mathrm{Ker%
}\Box $.

\begin{center}
\textbf{4. Three decomposition theorems }
\end{center}

In this paragraph we will consider the vector space Ker$\Box$ of all
infinitesimal harmonic transformations on a compact Riemannian manifold. The
following theorem is true.

\textbf{Theorem 4.1}. \textit{If the vector field }$\xi $\textit{\ is an
infinitesimal harmonic transformation on a compact Riemannian manifold }$%
(M,g)$\textit{\ then }$\xi $\textit{\ is decomposed in the form }$\xi =\xi
^{\prime }+\xi ^{\prime \prime }$\textit{\ where }$\xi ^{\prime }$\textit{\
is an infinitesimal isometric transformation and }$\xi ^{\prime \prime }$%
\textit{\ is an gradient infinitesimal harmonic transformation on }$(M,g)$%
\textit{. This decomposition is necessarily orthogonal with respect to the
global scalar product defined on }$(M,g).$

\textbf{Proof}. The vector space Ker$\Box \cap \mathrm{Ker\;d}^{\ast }$ of
all infinitesimal isometric transformations on a compact Riemannian manifold 
$(M,g)$ is a subspace of the finite-dimensional vector space Ker$\Box $ (see
Exp. 1). On the other hand it is well known (see [16], p. 205) that by
virtue of the Fredholm alternative vector spaces $\mathrm{Im\;d}$ and $%
\mathrm{Ker\;d}^{\ast }$ are orthogonal complement of each other with
respect to the global scalar product on compact Riemannian manifold $(M,g)$,
i.e. $\Omega ^{1}(M)=\mathrm{Ker\;d}^{\ast }\oplus \mathrm{Im\;d}$.
Therefore the vector space Ker$\Box \cap \mathrm{Ker\;d}$ of all
infinitesimal gradient harmonic transformations must be an orthogonal
complement of Ker$\Box \cap \mathrm{Ker\;d}^{\ast }$ with respect to the
whole space Ker$\Box $. This vector subspace consists of all gradient vector
fields $\nabla F$ such that $\nabla _{i}(\Delta F)=2R_{i}^{j}\nabla _{j}F$
for smooth scalar functions $F:M\rightarrow \mathbf{R}$.

\textbf{Remark}. The last result was known (see [23]) in the case of a
compact Einstein $n$-dimensional $(n\geq 2)$ manifold $(M,g)$ with constant
scalar curvature $s$.

Now we shell formulate the decomposition theorem of an arbitrary
infinitesimal harmonic transformation on a compact Kahlerian manifold.

\textbf{Theorem 4.2}. \textit{If }$\xi $\textit{\ is a holomorphic vector
field on a compact Kahlerian manifold }$(M,g,\mathrm{J})$\textit{\ then }$%
\xi $\textit{\ is decomposed in the form }$\xi =\xi ^{\prime }+\mathrm{J}\xi
^{\prime \prime }$\textit{\ where }$\xi ^{\prime }$\textit{\ and }$\xi
^{\prime \prime }$\textit{\ are both infinitesimal isometric
transformations. This decomposition is necessarily orthogonal with respect
to the global scalar product defined on }$(M,g,\mathrm{J}).$

\textbf{Proof}. On a compact Kahlerian manifold $(M,g,\mathrm{J})$, where as
well known $\nabla \mathrm{J}=0$, a vector field $\xi$ on a compact
Kahlerian manifold is an infinitesimal harmonic transformation if and only
if $\xi$ is a holomorphic vector field (see Exp. 5). Therefore, by virtue of
Theorem 4.1 we have the orthogonal decomposition $\xi=\xi^{\prime }+\mathrm{%
grad} F$ where $\xi^{\prime }$ is an infinitesimal isometric transformation
and grad $F$ is a holomorphic vector field for some smooth scalar function $%
F $ on $(M,g)$. On the other hand it is well known (see Theorem 6.8 of
Chapter IV in [24]) that $\mathrm{J}X$ is an infinitesimal isometric
transformation if a holomorphic vector field $X$ is closed. Therefore we can
state that $\xi=\xi^{\prime }+\mathrm{grad}F=\xi^{\prime }+\mathrm{J}%
\xi^{\prime \prime } $ where $\xi^{\prime \prime }$ is an infinitesimal
isometric transformation.

\textbf{Remark. }Lihnerowicz has proved the following theorem (see [13]): A
holomorphic vector field $\xi $ on a compact Kahlerian manifold $(M,g,%
\mathrm{J})$ with constant scalar curvature is decomposed in the form $\xi
=\xi ^{\prime }+\mathrm{J}\xi ^{\prime \prime }$ where $\xi ^{\prime }$ and $%
\xi ^{\prime \prime }$ are both infinitesimal isometric transformation.
Theorem 4.2 is a generalization of this theorem.

By virtue of the Fredholm alternative we shall prove the following theorem.

\textbf{Theorem 4.3}.\textbf{\ }\textit{On a compact Riemannian manifold }$%
(M,g)$\textit{\ of dimension }$(n\geq 2)$\textit{\ with positive Ricci
curvature an arbitrary infinitesimal conformal transformation }$\xi $\textit{%
\ has the form }$\xi =\xi ^{\prime }+\mathrm{grad}F$\textit{\ where }$\xi
^{\prime }$\textit{\ is an infinitesimal isometric transformation and F is a
some smooth scalar function on }$(M,g)$\textit{\ such that the vector field }%
$\mathrm{grad}$\textit{F is an infinitesimal conformal transformation.
Moreover, if }$L_{\mathrm{grad}F}s=0$\textit{\ then manifold }$(M,g)$\textit{%
\ is isometric to sphere }\textbf{S}$^{n\,}$\textit{in a Euclidian space }%
\textbf{R}$^{n+1}.$

\textbf{Proof}. The vector space of all infinitesimal conformal
transformations on a compact Riemannian manifold $(M,g)\ $is a
finite-dimensional vector space and the vector space of all infinitesimal
isometric trasformations is a subspace of this vector space. On the other
hand, there does not exist a nonzero harmonic vector field on a compact
Riemannian manifold with positive Ricci curvature (see Theorem 2.3 of
Chapter II in [22]), then $\mathrm{Im\;d\ =\mathrm{Ker\;d}}$. Therefore on a
compact Riemannian manifold (by virtue of the Fredholm alternative) for an
arbitrary infinitesimal conformal transformation $\xi $ the following
decomposition is true $\xi =\xi ^{\prime }+\mathrm{grad}F$ where $\xi
^{\prime }$ is an infinitesimal isometric transformation and $F$ is a some
smooth scalar function on $(M,g)$ such that the vector field $\xi ^{\prime
\prime }=\mathrm{grad}F$ with local coordinates $g^{ik}\nabla _{k}F$ is an
infinitesimal conformal transformation. Then by direct computation, we
obtain $L_{\xi }g=L_{\xi ^{\prime }}g+L_{\mathrm{grad}F}g=L_{\mathrm{grad}%
F}g=2\nabla \nabla F$ and $div\;\xi =-\Delta F$. As a result we receive the
following equality $L_{\mathrm{grad}F}g=2\nabla \nabla F=-(\Delta F)g$ from
which we can conclude that the vector field $\xi ^{\prime \prime }=\mathrm{%
grad}F$ is an infinitesimal conformal transformation also. It well known, if
a compact Riemannian manifold ($M,g)$ of dimension $n\geq 2$ admits a
nonconstant scalar function $F$ such that $\nabla \nabla F=n^{-1}(-\Delta
F)g $ then ($M$, $g)$ is conformal to a sphere \textbf{S}$^{n\,}$in
Euclidean space \textbf{R}$^{n\,+\,1}$ (see [11]). If in addition we suppose
that $L_{\mathrm{grad}F}s=0$ then ($M,g)$ must be isometric to a sphere 
\textbf{S}$^{n+1}$ (see [11]).

\textbf{Remark. }The vector space of all infinitesimal conformal
transformations on (\textbf{S}$^{n},\overline{g})$ splits as the direct sum
of the vector space of all infinitesimal isometric transformations and the
vector space of gradient vector fields of first speherical harmonics of \ (%
\textbf{S}$^{n},\overline{g})$. In particular, the vector space of all
infinitesimal conformal transformations on (\textbf{S}$^{2},\overline{g})$
has dimension equal to 6 and admits decomposition in the sum of two
subspaces (see [4]). Three of the dimensions arise from $\overline{\nabla }F$
where $F$ is a spherical harmonic. The other three dimensions come from the
infinitesimal isometric transformations for the standard metric $\overline{g}
$ on \textbf{S}$^{2}$. Therefore our decomposition the vector space of
infinitesimal conformal transformations on a compact Riemannian manifold is
an analog of above decomposition on a sphere (\textbf{S}$^{n},\overline{g})$.

\begin{center}
\textbf{5. Ricci solitons }
\end{center}

Let $(g,\xi,\lambda)$ be a \textit{Ricci soliton} on a smooth $n$%
-dimensional manifold $M$ (see Exp. 4), where $g$ is Riemannian metric and $%
\xi$ is a smooth vector field on $M$ such that the identity 
$$
-2Ric=L_\xi g+2\lambda g\eqno(5.1) 
$$
holds for some constant $\lambda$ (see [1], p. 22; [2], p. 353). Ricci
soliton is called \textit{steady}, if $\lambda = 0$, \textit{shrinking}, if $%
\lambda < 0$, and, finally, \textit{expanding}, if $\lambda > 0$.

In case $\xi=\mathrm{grad}F$ for some smooth function $F:M\rightarrow\mathbf{%
R}$ the equation can be rewritten as 
$$
-Ric=\nabla\nabla F+\lambda g.\eqno(5.2) 
$$
and $(g,\xi,\lambda)$ is called a \textit{gradient Ricci soliton} (see [1],
p. 22; [2], p. 353). Moreover, $(M,g)$ is called a \textit{trivial Ricci
soliton}, if $F=const$ and hence $(M,g)$ is an Einstein manifold.

{By Example 4 a vector field $\xi$ that makes a Riemannian metric $g$ into a
metric of a Ricci soliton is necessarily an infinitesimal harmonic
transformation. In addition, by the first decomposition theorem a harmonic
transformation $\xi$ on a compact Riemannian manifold $(M,g)$ has the form $%
\xi=\xi^{\prime }+\xi^{\prime \prime }$ where $\xi^{\prime }$ is an
infinitesimal isometric transformation and $\xi^{\prime \prime }$ is a
gradient infinitesimal harmonic transformation on $(M,g)$. By these
propositions we can rewrite the identity (5.1) as } 
\[
-2Ric=L_\xi g+2\lambda g=L_{\xi^{\prime }+\xi^{\prime \prime }}g+2\lambda
g=2\nabla\nabla F+2\lambda g 
\]
{where $\xi^{\prime \prime }=\mathrm{grad}F$ for some smooth scalar function 
$F$. Now we can formulate the following proposition. }

\textbf{Theorem 5.1.} \textit{The vector field }$\xi $\textit{\ of any Ricci
soliton }$(g,\xi ,\lambda )$ \textit{on a compact smooth manifold M \ has
the form }{$\xi =\xi ^{\prime }+\xi ^{\prime \prime }$}\textit{\ \ where }$%
\xi ^{\prime }$\textit{\ is an infinitesimal isometric transformation and }{$%
\xi ^{\prime \prime }$}\textit{\ is a gradient infinitesimal harmonic
transformation \ and \ therefore }$(g,\xi ,\lambda )$\textit{\ is a gradient
Ricci soliton.}

\textbf{Remark.} By means of Perelman work [15] and previous others, see
Hamilton [7] for dimension two and Ivey [10] for dimension 3 we know that an
every compact Ricci soliton is a gradient Ricci soliton. And hence the
Perelman-Hamilton-Ivey propositions is a corollary of our theorem about
infinitesimal harmonic transformations on a compact smooth manifold.

{Now we take the divergence of the Ricci tensor of $(g,\xi,\lambda)$. By
using the equation (5.2) we have $div \mathrm{Ric}=g^{ij}\nabla_iR_{jk}=%
\frac{1}{2}\nabla_ks=-g^{ij}\nabla_i \nabla_j\nabla_k F=\Delta(\nabla_kF)$.
Using this equation and the Schur's Lemma $div \mathrm{Ric}=2\mathrm{d}s$ we
get } 
$$
\mathrm{d}s=2\mathrm{d}(\Delta F).\eqno(5.3) 
$$
{Then by means of the equation (5.3) we have } 
\[
\langle\xi,\mathrm{d}s\rangle=2\langle\xi,\mathrm{d}(\Delta
F)\rangle=2\langle\mathrm{d}^*\xi,\Delta F\rangle=2\langle\mathrm{d^*d}%
F,\Delta F\rangle=2\langle\Delta F,\Delta F\rangle 
\]
{that is equivalent to $\int_M\xi(s)dV=\int_M(\Delta F)^2dv\geq 0$. By means
of this inequity we can formulate the following proposition (see [20])}

\textbf{Theorem 5.2.} \textit{If a shirking Ricci soliton $(g,\xi,\lambda)$
on compact smooth manifold $M$ satisfies the condition $L_\xi s\leq 0$ then
this soliton is trivial.}

\textbf{Remark}. It is well known that a compact steady or expanding Ricci
soliton $(g,\xi,\lambda)$ is a gradient soliton (see [15]) and on the other
hand a compact gradient steady or expanding Ricci soliton is a trivial
soliton (see [8]). On the other hand every shirking compact Ricci soliton
when $n > 3$ and the Weyl tensor is zero is trivial (see [5] and [25]). But
there is the open problem (see [5], p. 11): Are the special conditions in
dimension $n \geq 4$ assuring that a shirking compact Ricci soliton is
trivial? Our Theorem 5.2 may be is one of possible answers to this question.

We have proved (see [17]) that on a compact Riemannian manifold $(M,g)$ does
not exist an infinitesimal harmonic transformation $\xi $ such that $\mathrm{%
Ric}(\xi ,\xi )<0$. And in addition if $\mathrm{Ric(\xi ,\xi )\leq 0}$ for
an infinitesimal harmonic transformation $\xi $ such that $\xi \neq 0$ then $%
\xi $ is a parallel vector field.

On the other hand, the vector field $\xi$ that makes a Riemannian metric $g$
into a metric of Ricci soliton must be a null-vector field if $\nabla\xi=0$.
These two facts can be used to formulate the following assertion (see [20]).

\textbf{Corollary}. \textit{A Riemannian metric $g$ on a compact smooth
manifold $M$ can not be metric of a Ricci soliton }$(g,\xi,\lambda)$\textit{%
\ if }$\mathrm{Ric}(\xi,\xi)<0$\textit{. If }$\mathrm{Ric}(\xi,\xi)\leq0$%
\textit{\ then }$(g,\xi,\lambda)$\textit{\ is a trivial Ricci soliton. }

\end{document}